\documentclass{baustms}

\usepackage{booktabs}

\citesort

\theoremstyle{bams}
\newtheorem{thm}{Theorem}[section]

\newtheorem{cor}[thm]{Corollary}

\theoremstyle{bamsdefn}
\newtheorem{defn}[thm]{Definition}

\newcommand{\dia}{$\diamondsuit$ }
\newcommand{\diaa}{$\diamondsuit\!$}
\newcommand{\diaas}{$\diamondsuit$\nobreak\hskip -0.3 pt s }
\newcommand{\diab}{$\blacklozenge$ }

\newcommand{\diabs}{$\blacklozenge$\nobreak\hskip -0.3 pt s }
\newcommand{\whbox}{$\square$ }
\newcommand{\blbox}{$\blacksquare$ }
\newcommand{\whboxx}{$\square$\hskip -0.2 pt}
\newcommand{\blboxx}{$\blacksquare$\nobreak\hskip -0.1 pt}

\begin{document}
\runningtitle{Perfect Colourings of Isonemal Fabrics by Thick Striping}
\title{PERFECT COLOURINGS OF ISONEMAL FABRICS BY THICK STRIPING}
\cauthor 
\author[1]{R.S.D.~Thomas}
\address[1]{St John's College and Department of Mathematics, University of Manitoba,
Winnipeg, Manitoba  R3T 2N2  Canada.\email{thomas@cc.umanitoba.ca}}

\authorheadline{R.S.D.~Thomas}



\begin{abstract}
Perfect colouring of isonemal fabrics by thick striping of warp and weft and the closely related topic of isonemal prefabrics that fall apart are reconsidered and their relation further explored. The catalogue of isonemal prefabrics of genus V that fall apart is extended to order 20 with designs that can be used to weave cubes with colour symmetry as well as weaving symmetry.
\end{abstract}


\maketitle

\section{Introduction}
\noindent Except for a finite list of interesting exceptions, Richard Roth \cite{R1} classified isonemal periodic prefabric designs into 39 infinite species. Coarser than the five previously defined genera \cite{C1} are three more general and easily described classes. 
Species 1--10 have reflection or glide-reflection symmetries with parallel axes and no rotational symmetry, not even half-turns. 
Species 11--32 have reflection or glide-reflection symmetries with perpendicular axes, hence half-turns, but no quarter-turns. 
Species 33--39 have quarter-turn symmetries but no mirror or glide-reflection symmetries. 
This taxonomy has been refined slightly and used in \cite{P1,P2,P3}, to which reference needs to be made, to determine the feasible symmetry groups and hence isonemal prefabrics.
As Roth observes beginning his subsequent paper \cite{R2} on perfect colourings, `[r]ecent mathematical work on the theory of woven fabrics' begins with \cite{ST}, which remains the fundamental reference.
In \cite{R2} Roth determines which fabrics---actually prefabrics---can be perfectly coloured by striping warp and weft.
In the paper \cite{P4}, to which reference also needs to be made, I reconsider thin striping in terms of Roth's taxonomy as refined in \cite{P1,P2,P3} and consider further the related question of which isonemal prefabrics of even genus fall apart.
Here I pursue the topic of thick striping and questions about isonemal prefabrics of genus V that fall apart.
Striping is partly explained in Section 2.
In Section 3 designs of species 11--32 are considered and in Section 4 those of species 33--39.
In Section 5 the extension to order 20 begun in \cite{P4} of the catalogue \cite{JA} of isonemal prefabrics that fall apart is completed.
Section 6 considers the two-colouring of woven cubes.

\section{Striping}
\noindent We turn now to the matter of perfect colourings of the strands of a prefabric with two colours, the subject of \cite{JA}, Roth's later weaving paper \cite{R2}, and \cite{P4}.
If all of the symmetries of a prefabric with coloured strands are colour symmetries, then the choice of the strand colours is said to be {\it perfect} (or {\it symmetric}).
There are only two ways to colour a prefabric that can result in perfect colouring other than the {\em normal} (dark warps, pale wefts).
Warps and wefts can be striped, that is be pale and dark, either {\it thinly}, that is alternately, or {\it thickly}, that is alternating in pairs: pale, pale, dark, dark, pale, pale.
Adapting a device from \cite{R2}, the colouring of a prefabric can be represented by seeming to extend strands outside the pattern to indicate which strands are pale or dark.
I adopt this convention as long as it is not completely obvious which strands are which (Figures 1, 2, and 13).

Striping warp and weft creates a checkerboard of cells that may be called redundant and irredundant, where the {\it redundant} are those where the same colour meets itself and the {\it irredundant} are those where the pattern colour is determined by the design.
In this language, the irredundant cells on the obverse have the complement of the design in predominantly dark rows and the colour of the design in predominantly pale rows, in both cases the reverse has the complementary colour.
\begin{figure}
\centering
\includegraphics{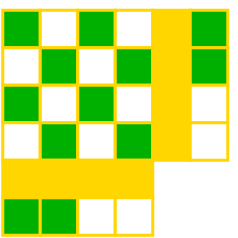}\hskip 10 pt
\raisebox{22 pt}{\includegraphics{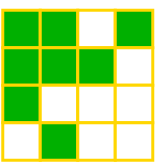}\hskip 10 pt
\includegraphics{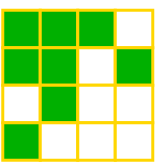}}

\hskip 10pt (a)\hskip 54 pt (b)\hskip 42 pt (c)
\caption{Plain weave.\hskip 10 pt a. One colouring by thick striping. \hskip 10 pt b. Obverse view. \hskip 10 pt c. Reverse view (i.e., other side as viewed in a mirror).}\label{fig1:}
\end{figure}
\begin{figure}
\centering
\includegraphics{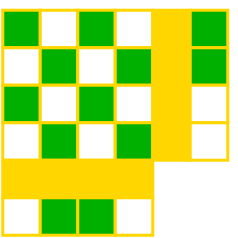}\hskip 10 pt
\raisebox{22 pt}{\includegraphics{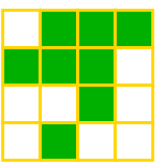}\hskip 10 pt
\includegraphics{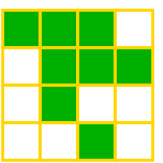}}

\hskip 10pt (a)\hskip 54 pt (b)\hskip 42 pt (c)
\caption{Plain weave.\hskip 10 pt a. Second colouring by thick striping. \hskip 10 pt b. Obverse view. \hskip 10 pt c. Reverse view.}\label{fig2:}
\end{figure}
Consider the effect of striping the warp and weft of plain weave thickly.
The checkerboards of redundant cells in both of the distinct stripings of Figures 1 and 2 are that of plain weave doubled (box weave 4-3-1). 
The resulting pattern in both cases is the simplest houndstooth 4-1-2*, a prefabric that falls apart if coloured normally but, as is illustrated here, is weavable as a fabric with striped strands.
This striping illustrates that the patterns obtained by striping warp and weft of isonemal fabrics are the designs of prefabrics that fall apart \cite[Lemma 3]{JA}.
In consequence of this fact together with the known conditions for isonemal fabrics that fall apart (\cite{WD}, \cite{CRJC}), if the pattern of a fabric obtained by striping warp and weft is the design of an isonemal prefabric, then the prefabric is of genus II, IV, or V with no overlap with genus I or III and with one quarter of the cells in half the rows dark.
The pattern arising from the striping of warp and weft of an isonemal fabric does not need to {\it be} the design of an isonemal prefabric \cite[\S 3]{P4}.

If there is to be any hope of perfect colouring, then the colouring, stri\-ping, must be chosen so that the colour symmetries of the fabrics map redundant cells to redundant cells and irredundant cells to irredundant cells, or as Roth puts it `preserve' them.
Which half of the cells are to be redundant and which half are to be irredundant is a choice to be made.
There are therefore two ways to stripe the same fabric thinly and two thickly, taking reversal of colours not to be a different striping.
Preserving the two classes of cell rules out as symmetries, for any striping, glide-reflections with axes not in mirror position, translations $(x, y)$ with $x$ and $y$ not integers of the same parity, and half-turns with centres not at the centre or corner of a cell but not half-turns with centres in those two positions, not other translations and not mirror symmetries.
It rules out, for thin striping, quarter-turns with centres not at the centre of a cell.
It rules out, for thick striping, translations $(x, y)$ with $x$ and $y$ odd or with $x+y$ not divisible by 4.
The two theorems on what can be thinly striped, essentially proved by Roth, although not stated in terms of his taxonomy, are (from \cite{P4}) as follows.
\begin{thm}
Isonemal periodic prefabrics of order greater than $4$ and of every species with symmetry axes that is not ruled out by the placement of glide-reflections can be perfectly coloured by thin striping: $1_m, 2_m, 3, 5_o, 5_e, 6, 7_o, 7_e, 8_e$, $9, 11, 13, 15, 17, 18_s, 19, 21$--$23, 25$--$27, 28_e, 28_n$, and $29$--$31.$
\end{thm}

\begin{thm}
Isonemal periodic prefabrics of order greater than $4$ with quarter-turn symmetry can be perfectly coloured by thin striping if and only if they are of species $36_s$.
\end{thm}

For thick striping, the translation constraint eliminates as candidates all twills and twillins (genera I and II), having as symmetries, as they do, translations from strand to adjacent strand ($x = 1$ or $y=1$), taking with them species 1--10, 12, 14, 16, 18, 20, 23, 24, 26, 28, 31, and 32 from the first 32 species.
Order, being the minimal $p$ such that $(p, 0)$ is a translation, must be divisible by 4 eliminating, in addition to those already forbidden, only species $27_o$.
In addition, for thick striping, quarter-turn centres must be at cell corners and fall in the centre of blocks of four cells (redundant or irredundant), and half-turn centres must also be at cell corners but can fall in the centre or at the corners of blocks of four cells (redundant or irredundant).
Roth has shown \cite{R2} that these modest necessary conditions are also sufficient to allow the two sorts of striping.

\section{Symmetry with Perpendicular Axes}
\noindent To use thick striping to produce perfect colouring, symmetry groups with glide-reflections must have axes in mirror position, translations $(x, y)$ with $x$ and $y$ even and with $x+y$ divisible by 4, any quarter-turn centres at cell corners that fall in the centre of blocks of four redundant or irredundant cells, and half-turn centres at cell corners that fall in the centre or at the corners of blocks of four redundant cells.
As mentioned, these constraints eliminate species 1--10 and 12, 14, 16, 18, 20, 23, 24, 26, $27_o$, 28, 31, and 32 from those with perpendicular axes.
Half-turns with centres not at cell corners are ruled out, and those at cell corners have to be spaced to fall at the centre or corner of the blocks of four redundant and irredundant cells.
A mirror must pass alternately through centres and corners of such blocks.
Other mirror positions are ruled out, again a matter of spacing.
Finally, glide-reflections are ruled out unless they lie in the position just described for mirrors, and their glides must be an even multiple of $\delta$ in order to preserve the doubled checkerboard of blocks of redundant cells.
We need to examine species 11, 13, 15, 17, 19, 21, 22, 25, $27_e$, 29, and 30 to see which might produce perfect colouring by having strands striped thickly.

Roth points out that when $G_1$ is of type $pgg$, the edges of its lattice unit must be allowable translations, that is with components divisible by 4.
But from \cite{P2} they are $2a$ and $2b$ with $(a, b)=1$ and so cannot be divisible by 4.
This rules out the species from 11 to 16.

The group type $pmg$ is not so restrictive.
If the quarter $G_1$ lattice unit is odd by odd in $\delta$ units, then the lattice-unit edges are not allowable translations, but if it is odd by even and the even direction is perpendicular to the mirrors, that is, it is the distance between the mirrors, then double the glide will be the necessary multiple of $\delta$ divisible by 4.
This bans species species $17_o$ and $19_o$ but allows $17_e$ and $19_e$, provided parameters $a$ and $b$ are used only the right way around and not with the other orientation.
The positive possibility is illustrated by 8-19-7 (Roth's $17_e$ example \cite[Figure 9a]{P2}) and the negative by 8-19-4 ($19_e$ example \cite[Figure 9c]{P2}), which has mirrors $\delta$ apart.
Satisfaction of the same restriction is forced on all designs of species 21, for example 8-7-2 \cite[Figure 10]{P2}.
Species 22, on the other hand, is forced to violate the restriction, the distance between mirrors being always odd in $\delta$.
That concludes the groups of type $pmg$, allowing $17_e$, $19_e$, and 21.

The only species of crystallographic type $pmm$ that might work is 25.
All such designs have translations (edges of the lattice unit) of even length in $\delta$ because of the standard isonemal spacing of the quarter lattice unit.
Accordingly all such fabrics can be thickly striped, which is not to say that every thick striping is a perfect colouring.
The striping has to be done to be compatible with the location of the \diaas\hskip -3pt .

The remaining possibilities are of crystallographic type $cmm$.
As for thin striping, stripability depends on the translations $(M, N)$ and $(N, M)$, that are the sides of the rhombic lattice unit.
Whereas for thin striping their components needed just the same parity, for thick stiping they must be even and $M+N$ and $M-N$, the length and width of the rhombs in $\delta$ units, must be divisible by 4. 
So only the even-even spacing of species $27_e$ (example \cite[Figure 12b]{P2}), also used in species 29 (Roth example 16-2499 \cite[Figure 6]{R1} = \cite[Figure 15a]{P2}), is acceptable.

These facts can be summarized in a theorem containing no information not in \cite{R2}, but \cite{R2} is not written in terms of Roth's symmetry-group types.
\begin{thm}
The species of isonemal prefabrics that can be perfectly coloured by thick striping of warp and weft are both $17_e$ and $19_e$ provided that the distance between neighbouring mirrors is even in $\delta$ units, $21, 25, 27_e$, and $29$. 
\end{thm}
\noindent The species that allow thick striping are a proper subset of those that allow thin striping (Cf. Theorem 2.1).
\begin{cor}
Isonemal prefabrics that can be perfectly coloured by thick striping of warp and weft have orders divisible by eight.
\end{cor}
\noindent This corollary is a direct consequence of the constraints on orders of the species \cite[\S 12]{P2} mentioned in the theorem. In fact, prefabrics of species $27_e$ and $29$ have orders divisible by 16.

More interesting and less obvious is the correspondence between prefabrics of the species that allow thick striping, $17_e$, 21, 25, $27_e$, and 29, and some prefabrics of the species $18_s$, 23, 26, $28_e$ and $28_n$, and 31, which are the only species that can both be doubled and remain isonemal \cite[Theorem 4]{P2} and be perfectly coloured by thin striping \cite[Theorem 2.1]{P4}.
Such prefabrics can be doubled because they have centres of half turns on their mid-lines.
When they are doubled they become $17_e$, 21, 25, $27_e$, and 29 respectively.
When one is striped thinly {\it and} doubled, what results is a thick striping.
Since the existence of some prefabrics that can be thickly striped shows that the types they fall in can be thickly striped, we have an interesting direct way to see this list.
But $19_e$ is missing from the list. 
It would be the result of `doubling' a non-existent stripable species $20_s^*$ with the spacing of $18_s$ but with \diab on mid-lines of strands where they are impossible (as was mentioned when species 19 and 20 were discussed in \cite{P2}).
After the `doubling', when \diab {\it can} be at the cell-corner position formerly the centre of a cell, we have fabrics of type $19_e$.

In view of the non-existence of a species $20_s^*$, one can wonder whether fabrics of species $19_e$ really can be thickly striped, but fabric 8-11-4 of species $19_e$ becomes 8-5-3* when thickly striped.

When fabrics are thinly striped, there is a marked tendency to stripiness of the patterns produced, so strong a tendency that designs of isonemal prefabrics are the exception rather than the rule.
Since the only species with perpendicular axes that we know will produce an isonemal design when striped thinly are 11, 22, and 30, we have reason to expect non-isonemal designs to be produced by thick striping, especially as so many thick stripings are doubled thin stripings, and that expectation is not disappointed.
Prefabrics of species $17_e$, 25, and $27_e$ have $H_1$ of crystallographic type $p2$ with no warp-to-weft transformation. 
One can expect such patterns to have stripes.
Examples show that prefabrics of the other types, $19_e$, 21, and 29, when thickly striped can fail to be designs of isonemal prefabrics because the pairs of strands of  the same colour cannot be interchanged or because the pale and dark pairs of parallel strands cannot be interchanged although there are warp-to-weft transformations.
These examples show that there can be no theorem like Theorem 3 of \cite{P4} assuring us that whole species of fabrics with perpendicular axes produce designs of isonemal prefabrics when thickly striped.

We can say the following.

\begin{thm}
The design of an isonemal prefabric of order greater than $4$ and of pure genus V, which includes all those that look thickly striped, is of species $21, 29,$ or $37$. 
\end{thm}
\begin{proof} These species were determined by Roth \cite{R1} to be the only non-exceptional prefabrics of pure genus V.
\end{proof}
\begin{cor}
If an isonemal fabric perfectly coloured by thick striping is the design of an isonemal prefabric, then the prefabric is of species $21, 29$, or $37$.
\end{cor}
\begin{proof} Because of the alternation in pairs of rows 3/4 dark and 3/4 pale, such a prefabric must be of genus V and no other.
\end{proof}
\begin{cor}
The order of the design of an isonemal prefabric with symmetry axes, of order greater than $4$, and of genus V is a multiple of $8$.
\end{cor}
\begin{proof} Orders of species 21 and 29 are divisible by 8 and 16 respectively \cite{P2}.
\end{proof}

Corollary 3.5 accounts for the known lack of such prefabrics of order 12 and guarantees such a lack at order 20.
Prefabrics of species 37 with quarter-turn symmetry {\em can} have order divisible by 4 \cite{P3}; we return to species 37 in the next section.
Prefabrics of species 21 and 29 have nothing to contribute to the extension to order 20 \cite{P4} of the catalogue of prefabrics that fall apart \cite{JA}.

It can be noted that 4-1-1*, the only isonemal prefabric of order 4 that falls apart and whose design is produced by thin striping, is of species $23_o$.
It can be doubled to produce 8-3-1* of species 21, and the only other thickly striped prefabric of order 8, 8-9-1* is a variation of 8-3-1* in the same species.
Similarly 8-5-1* of species $23_e$ doubles to 16-51-1* of species 21.
Variations of 16-51-1* in its order and species are 153-1*, 291-1*, 291-2*, 531-2*, 531-3*, and 2193-1*.
Likewise, 8-5-3* of species 31 doubles to 16-51-2* of species 29, which has variations there 153-2*, 291-3*, 291-4*, 531-1*, 531-4*, 561-1* and 2193-2*.
Only the two prefabrics 16-51-1* and 16-51-2* can be produced directly by doubling.

The above prefabric designs can all be produced by thickly striping isonemal fabrics.
This is not, however, generally possible.
The key fact that allows this to be seen is that all of the species with symmetry axes that allow thick striping determined in Theorem 3.3 have mirror symmetry.
Consider the isonemal prefabric of order 32 and species 21 that falls apart and is illustrated in Figure 3a, a variation of 16-85-1* doubled.
\begin{figure}
\centering
\includegraphics{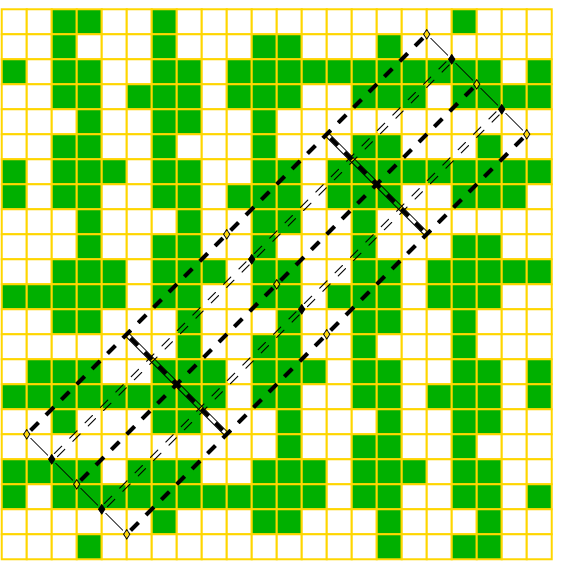}\hskip 10 pt
\includegraphics{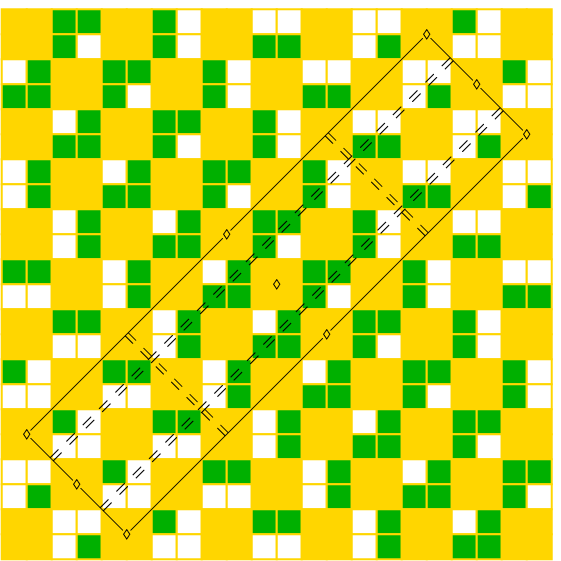}

(a)\hskip 2.1 in (b)
\caption{a. A prefabric of species 21 and order 32 that falls apart.\hskip 10 pt
b. Partial construction of a fabric (redundant cells neutral) that could be coloured to look like the design of (a).}\label{fig3:}
\end{figure}
If one wants to go from this prefabric that falls apart to the design of a fabric that could be coloured by thick striping to look like it, one will have to reverse the colours of the irredundant cells in the predominantly dark rows and not those in the predominantly pale rows so that the striping would give the starting design.
But when this is done with the prefabric of Figure 3a, the partial design that results, shown in Figure 3b, has nowhere that a mirror would fit among the irredundant cells determined, no matter what was done with the formerly redundant cells.%
\footnote{This example, like the one in \cite[\S 5]{P4}, disproves \cite[Theorem 2]{JA}.}
The situation is quite different in the next section.

\section{Quarter-turn Symmetry}
\noindent For thick strand striping to produce perfect colouring of prefabrics with quarter-turn symmetry, the components of the sides of the lattice units of level $i$, $M_i$ and $N_i$, must both be even and $M_i + N_i$ divisible by four.
The first condition requires level three or four, and the second, since $M_3 + N_3 = 2(M_1 + N_1)$, which is twice an odd number, requires level four, where 
$M_4 + N_4 = 2(M_2 + N_2)$, twice the sum of even numbers.
We consider the order-20 designs of Figures 4a, 4b, and 5 of the level-4 species $33_4$, $35_4$, and 37 respectively.
\begin{figure}
\centering
\includegraphics{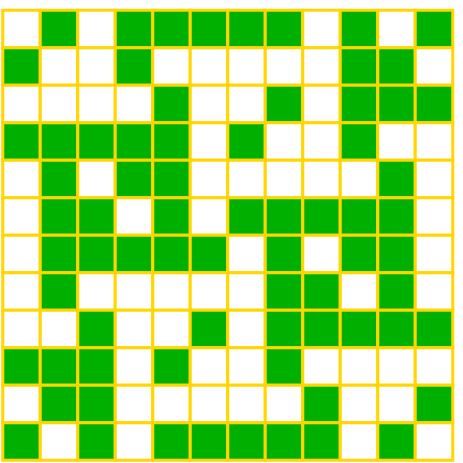}\hskip 10 pt
\includegraphics{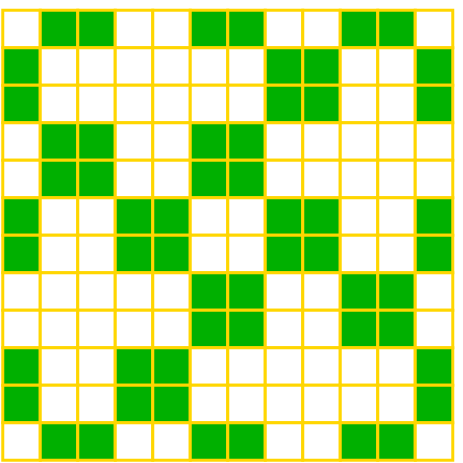}

(a)\hskip 128 pt (b)
\caption{a. Order-20 example of species $33_4$ \cite[Figure 13b]{P3}.\hskip 10 pt b. 10-85-1 doubled of order 20 and species $35_4$ \cite[Figure 13a]{P3}.}\label{fig4:}
\end{figure}
\begin{figure}
\centering
\includegraphics{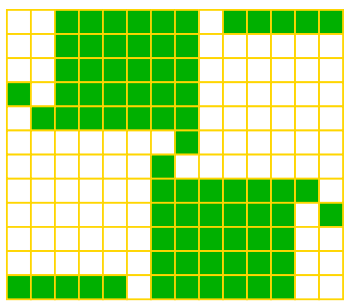}
\caption{Order-20 example of type 37 \cite[Figure 15]{P3}.}\label{fig5:}
\end{figure}

It must be kept in mind that, because prefabrics with order greater than 4 and quarter-turn symmetry lack all reflective symmetry, the mirror images of all prefabrics in the remainder of this paper are distinct from those illustrated, but they are systematically ignored.

The striping must be chosen so that the centres of quarter-turns are in the centre of blocks of four redundant or irredundant cells.
These centres of quarter-turns form two lattices, one of lattice-unit corners, the other of lattice-unit centres in what I have called chain-mail relation (\cite[Figure 5b]{P3}
or \cite[Figure 1d]{P4}).
The logical possibilities are shown in Table 1.
Because the movement from the corner of a level-four lattice to its centre is twice an even number of cell widths in one vertical or horizontal direction and then twice an odd number of cell widths in a perpendicular direction, $2(M_1 + N_1)= M_3 + N_3$, the two lattices of centre types specified in Table 1 fall one in the centres of redundant blocks and the other in the centres of irredundant blocks.
Because half the lattice-unit side is the hypotenuse of a right triangle with odd other sides, $M_2$ and $N_2$, the half-turns with centres at mid-side fall at the corner rather than at the centre of a block of redundant cells; this fact is needed for the proof of Theorem 4.1.
In lines 1 and 2 of Table 1 giving species $33_4$ and $35_4$, it does not matter which lattice is which, but for species 37, given by lines 3 and 4, the results in general differ.
Let lattice 1 fall in redundant blocks and lattice 2 fall in irredundant blocks.
Figures 6%
\footnote{The new catalogue numbers in the caption will be discussed in Section 5.}
and 7 illustrate species-$33_4$ and species-$35_4$ fabrics perfectly coloured by thick striping.
Which strands are pale and dark is sufficiently evident that it is not indicated separately in the figures.
Figure 7 illustrates that the phenomenon of the disappearing pattern (cf. \cite[Figure 3]{P4}) can occur with rotational symmetry too.
The symmetry groups displayed in Figures 6a, 7a, 8a, and 9a are those of the underlying fabric, not of the diagram itself.
The symmetry groups displayed in Figures 6b, 8b, and 9b are those of standardly coloured isonemal prefabrics with the appearance of those figures since conventions for such display are readily available.
The patterns are being treated as designs.
No symmetry group is marked on Figure 7b since it is so obviously not the design of an isonemal prefabric.
\begin{table}
\caption{Rotation assignments for two level-four lattices.}\label{tab:1}
\centering
\begin{tabular}{cccc}
\toprule
\multicolumn{1}{c}{Lattice 1} & Lattice 2 \\
\multicolumn{1}{c}{(redundant} & (irredundant & {Species}\\
\multicolumn{1}{c}{blocks)} &  blocks)\\
\midrule
\whbox & \whbox & $33_4$  \\
\blbox &\blbox  & $35_4$  \\
\whbox & \blbox & 37  \\
\blbox & \whbox & 37  \\
\bottomrule
\end{tabular}
\end{table}
\begin{figure}
\centering
\includegraphics{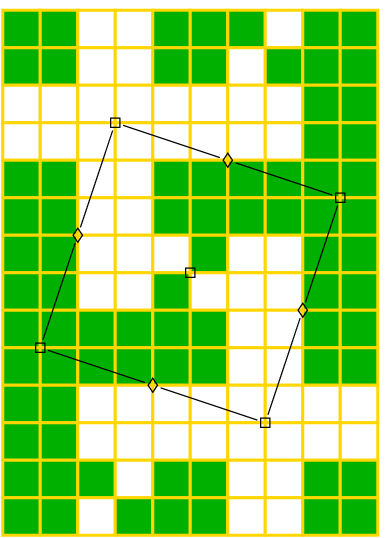}\hskip 10 pt
\includegraphics{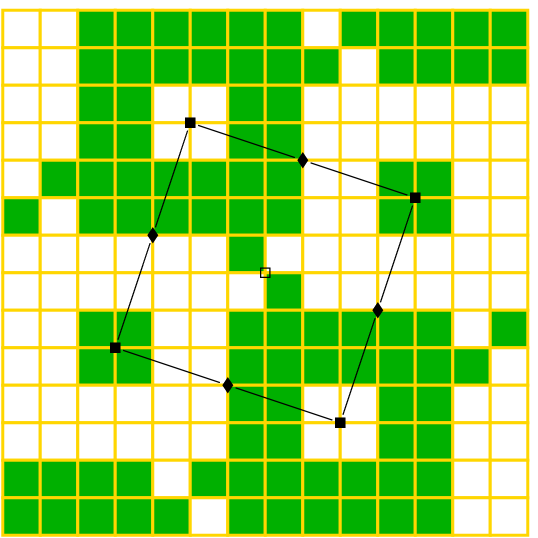}

(a)\hskip 130 pt (b)\phantom{xxxx}
\caption{Patterns of colouring by thick striping of the species-$33_4$ example of Figure 4a.\hskip 10 pt 
a. Obverse (20-787-2*)
with $G_1$ of the underlying fabric marked.\hskip 10 pt b. Reverse (20-4147-2*) with $G_1$ of both patterns marked.}\label{fig6:}
\end{figure}
\begin{figure}
\centering
\includegraphics{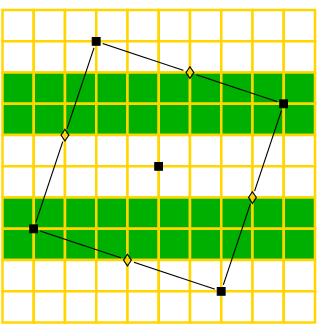}\hskip 10 pt
\includegraphics{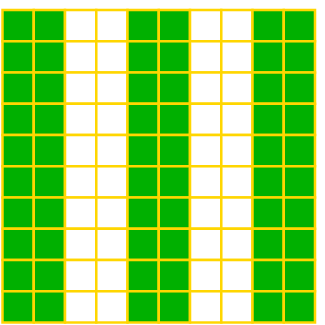}

(a)\hskip 90 pt (b)
\caption{Patterns of colouring by thick striping of the species-$35_4$ example of Figure 4b.\hskip 10 pt a. Obverse with $G_1$ of the underlying fabric marked.\hskip 10 pt b. Reverse.}\label{fig7:}
\end{figure}

Figures 6 and 7 also illustrate that the result of thick stiping may be and may not be the design of an {\it isonemal} prefabric that falls apart.
But we still have two lines of Table 1 to deal with.
Applying to the species-37 fabric of Figure 5 the scheme of line 3 of Table 1 gives us Figure 8, and the scheme of line 4 gives us Figure 9 as obverse and reverse patterns of the fabric.
$G_1$ of the underlying fabric is marked on each obverse, and $G_1$ of a design with the appearance of the patterns is marked on each reverse.
\begin{figure}
\centering
\includegraphics{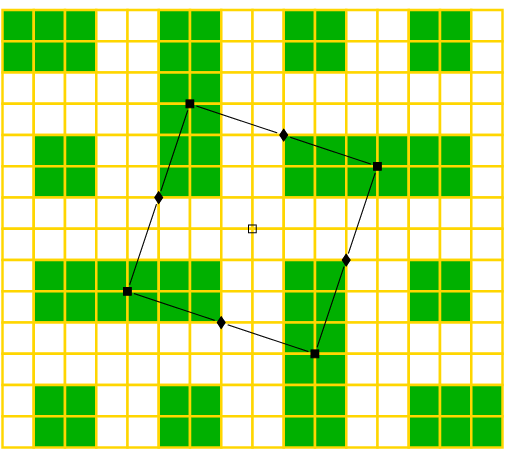}\hskip 10pt \includegraphics{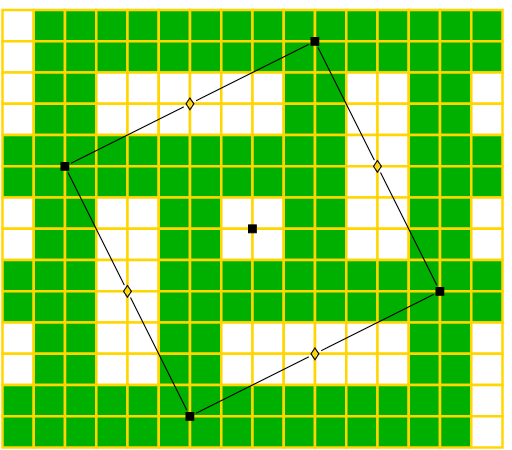}

(a)\hskip 142 pt (b)
\caption{A colouring by thick striping of the species-37 example of Figure 5. \hskip 10pt a. Obverse with $G_1$ of the underlying fabric marked.\hskip 10pt b. Reverse with $G_1$ of both patterns marked.}\label{fig8:}
\end{figure}

\begin{figure}
\centering
\includegraphics{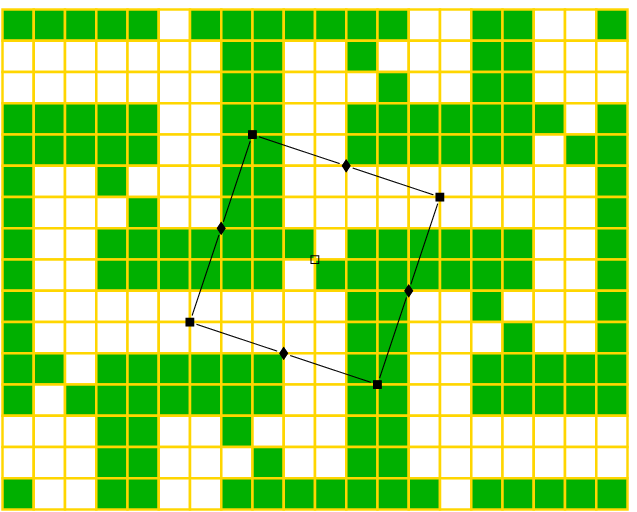}\hskip 5 pt\includegraphics{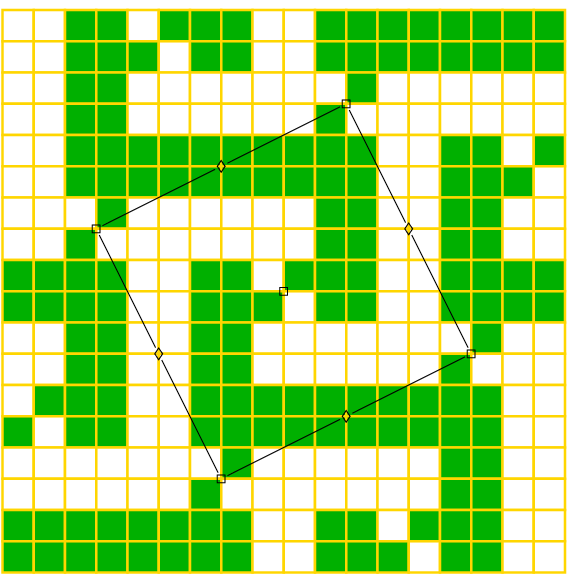}

\hskip 10pt(a)\hskip 169pt (b)
\caption{Second colouring by thick striping of the species-37 example of Figure 5. a. Obverse with $G_1$ of the underlying fabric marked. b. Reverse with $G_1$ of both patterns marked.}\label{fig9:}
\end{figure}
\begin{thm}{If the pattern of a perfect colouring of an isonemal fabric of order greater than $4$ with quarter-turn symmetry obtained by striping its strands thickly is the design of an isonemal prefabric, then the symmetry group of the prefabric is of type $37$.}
\end{thm} 
\begin{proof} The isonemal fabric coloured must be of level four.
We can reuse the lattice specifications of Table 1 for quarter-turn centres in the centres of blocks of 4 cells of the coloured pattern.
In lattice 1 \whbox can be seen to be impossible by considering the cells of the block itself or the four redundant blocks at its corners; in the design of an isonemal prefabric, \whbox cannot fall in the centre of four cells of the same colour.
This eliminates lines 1 and 3 of Table 1 as possibilities.
In lattice 2 \blbox can be seen to be impossible by considering the redundant blocks surrounding those irredundant blocks in which the \blbox might be; those related by such half turns are the same colour but the two pairs are opposite in colour.
In the design of an isonemal prefabric, \blbox cannot relate cells of different colours.
This eliminates lines 2 and 3 of Table 1 as possibilities.
Only line 4 remains, and its configuration of quarter-turn centres at level 4 characterizes Roth type 37.
\end{proof}

An alternative proof would use Theorem 3.3.
The theorem shows that the patterns of Figures 8 and 9 cannot be designs of isonemal prefabrics.
The designs, that is, patterns viewed as designs, of Figures 8 and 9 have symmetry groups of appropriate crystallographic types but at level 5 (cf.~\cite[Lemma 9]{P3}) and so not transitive on strands; the way centres of half-turns and of quarter-turns line up vertically and horizontally in Figures 8b and 9b indicates the trouble.
Some adjacent pairs of strands are interchangeable by them but are not related to the next adjacent pairs.

More constructive considerations of what happens in the thick striping illustrated in Figures 6--9 are these.

\noindent 1. \blbox is converted to \dia by the striping whether in a redundant or an irredundant block. 
If a redundant block surrounds \blbox then the irredundant blocks around it are no longer related by \blboxx , just by \diaa , because either those above and below or those beside it are complemented by the striping.
If an irredundant block surrounds \blboxx , then \blbox is preserved in that block, but the striping forces \whbox on the redundant blocks above, below, and beside it; this contradiction destroys \blbox but preserves  $\hbox{\blboxx\hskip 0.7 pt}^2 = \hbox{\whboxx\hskip 0.7 pt}^2 = \hbox{\diaa}$.

\noindent 2. \whbox is invariant in irredundant blocks and is converted to \blbox by the imposition of a redundant block.
Making a block irredundant leaves the irredundant blocks in the rows and the columns in which it lies unchanged or complements both leaving them related by \whbox if they were before.
Making a block redundant complements either the irredundant blocks in the rows or the columns in which it lies, changing what was related by \whbox in these rows and columns to being related by \blboxx .
That these two effects extend to the whole plane depends on the position of the centre of quarter-turn within the checkerboard of redundant blocks and on the effects of the complementation on the irredundant blocks.

\noindent In Figure 7, the symmetry group of the pattern becomes of type $p2$, which the symmetry group of an isonemal prefabric cannot be.
There is nothing to relate wefts to warps, and that is what is visible.
The explanation is general to the species and also covers the thin striping of species $36_s$, whose stripiness when striped \cite[Figures 7, 8]{P4} is now explained.
Rotation with side-reversal is still a symmetry of the {\it fabric,} but side-reversal means just that.
The stripes of Figure 7a are rotated and appear (reflected by the convention on display of reverse sides \cite[\S 1]{P1}) on the other side in Figure 7b, and the motifs of Figure 8a are rotated and appear on the other side in Figure 8b, likewise Figures 9a and 9b.
And vice versa in all cases.
The reversed colour of the motifs in Figures 7 and 8 and non-reversal in Figure 9 reflects the location of \blbox in irredundant blocks in 7 and 8 and in reduntant blocks in 9.
The effect of the change of \blbox to \dia is different for fabrics of species 37.
The quarter-turns \blbox become \diaas at mid-sides of a new larger lattice unit---in Figure 9b, just as though side-preservation were imposed.
As far as each side separately is concerned, that is, for design purposes, the half-turns at mid-side of the former lattice unit simply disappear.

\begin{thm}{The pattern of a fabric of order greater than $4$ with quarter-turn symmetry and perfectly coloured by thickly striping strands is the design of an isonemal prefabric that falls apart if and only if the fabric is of species $33_4$.}
\end{thm} 
\begin{proof}
Only if. Since an isonemal prefabric is known from the proof of the previous theorem to be of species 37 with \whbox in the centre of irredundant blocks and \blbox in the centre of redundant blocks, it suffices to see where they come from. Cf. Figure 6.
A \whbox in the centre of an irredundant block comes from a \whbox in the centre of a block of four cells in the fabric.
A \blbox in the centre of a redundant block might come from a \blbox in the center of a block or a \whbox in the centre of a block, since imposition of the redundancy of the block forces \whbox to be a \blboxx .
Observation 1 of the previous paragraph shows that if it had been \blboxx , it would have been destroyed.
Accordingly a \blbox does not come from a \blboxx .

No contradiction arises from the imposition of redundancy on a block of cells surrounding a \whboxx .
The four irredundant blocks surrounding its block are images of one block but with colouring the same in opposite pairs, one pair the colour-complement of the other pair in accordance with \whbox (no $\tau$).
When the complementation occurs on account of the striping, they cease to be colour-complements and become all the same in accordance with the \blbox in the centre of the redundant block.

Since both the \blbox and \whbox symmetries of the prefabric come from \whbox symmetries of the fabric, the two level-4 lattices in the fabric are those characteristic of species $33_4$.
There could perhaps be more symmetries.
There are two ways in which a type-$33_4$ group present might be a proper subgroup of the symmetry group.
One is that the subgroup has side components $M^{\prime}_4$, $N^{\prime}_4$, that are multiples of the components of another type-$33_4$ group: $(M^{\prime}_4, N^{\prime}_4) = (pM_4, pN_4)$ relating level-4 lattice units based on, say, $M^{\prime 2}_1 + N^{\prime 2}_1 = p^2q^2$ to those with $M^2_1 + N^2_1 = q^2$.
But such multiplication would spoil the required relative primality of $M^{\prime}_1$ and $N^{\prime}_1$.
This does not happen.
The other way is that the prefabric has the specified level-4 lattice unit but the underlying fabric has more symmetry, being of species $33_3$, $34$, or 39.
But fabrics of these species cannot be {\it perfectly} coloured by thick striping.
This too does not happen.
There cannot then be more symmetries.
The fabric is of species $33_4$.

If. Species $33_4$ has two lattices of \whbox in chain-mail relation.
Either can be taken as centres of redundant blocks, the other as centres of irredundant blocks.
In the irredundant blocks, \whboxx s remain \whboxx s, and in a redundant block a \whbox is converted to \blboxx . The configuration of quarter-turn centres characteristic of species 37 results.
As remarked in the third paragraph of this section, the \diaas at mid-sides of the lattice units fall on corners of redundant and irredundant blocks, and so on the edge of the thick stripes.
The striping accordingly turns them into \diabs as required for species 37.
\end{proof}

The theorem would be false without the restriction to perfect colouring.
Fabrics of species $33_3$, 34, and 39 can be thickly striped to produce patterns of isonemal prefabrics of species 37 that fall apart, but the colourings are not perfect.
Examples are, respectively 10-55-2 \cite[Figure 11b]{P3}, 10-107-1 \cite[Figure 8a]{P3}, and 10-93-1 \cite[Figure 7]{P3}.
These fabrics have \whboxx s in the right places to turn into the \whboxx s, \blboxx s, and \diabs required for a group of type 37, but some of the symmetries of the fabric vanish altogether from the coloured pattern: \dia from 10-55-2 of level 3 and from 10-107-1 of level 2 and both \blbox and \diab (and some \whboxx s) from 10-93-1 of level 1.
Examples of resulting coloured patterns are respectively what will shortly be catalogued as prefabrics 20-4147-2*, 20-4489-2*, and 20-4371-1*.

Because thin strand stripings of species $36_s$, if doubled, become thick strand stripings of fabrics of species $35_4$, the fact that a fabric of species $35_4$, when thickly striped, cannot have a pattern that is the design of an isonemal prefabric shows the same to be true of fabrics of species $36_s$ when thinly striped.
This is a different proof of \cite[Theorem 4]{P4}.

\section{Pre-fabrics that fall apart}
\noindent While it may be surprising that it has not previously been noted that all patterns resulting from thick striping of fabrics of species $33_4$ are isonemal prefabrics that fall apart, the reason is that no {\it example} of an isonemal prefabric with rotational symmetry that falls apart has been noted prior to Figure 6 on account of its comparatively large order, 20, except for 4-1-2*, discounted because it is so obviously exceptional.
The catalogue of isonemal prefabrics that fall apart in \cite{JA} extends only to order 16.
The designs of species 37 of order 20 are not so plentiful that it is infeasible to illustrate all 30, which I do in Figures 10--12, except for the two in Figure 6.
Some explanation needs to be given of how it is possible to produce such an exhaustive list with none of the trial and error characteristic of previous such figures \cite{C1, C2, WD, JA}.
Each cell has an orbit under the symmetry group.
In these designs, the orbits of only five arbitrarily colourable cells cover the plane, so that all of the colour choices can easily be investigated.
When this is done, the 30 patterns of Figures 6 and 10--12 result, together with two versions of the houndstooth 4-1-2*, i.e., the pattern of Figure 1b 
and its reverse, which is a mirror image of Figure 1b. 
These two patterns (each the other's reverse) will occur in any such production because every group of type 37 is a subgroup of the symmetry group of 4-1-2*.
\begin{figure}
\centering
\includegraphics{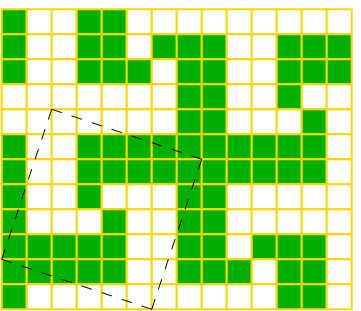}\hskip 10pt \includegraphics{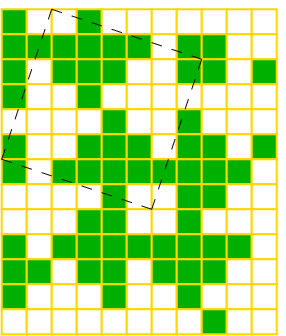}\hskip 10pt \includegraphics{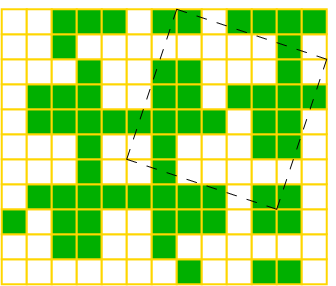}

(a) 787-1*/4147-1* \hskip 23pt (b) 2329-1*/4249-1* \hskip 13pt (c) 2329-2*/4249-2*

\vspace {4 pt}
\includegraphics{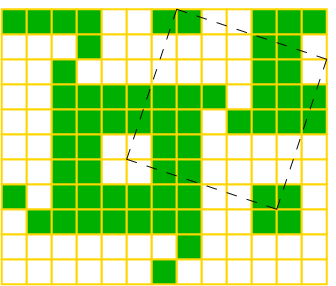}\hskip 10pt \includegraphics{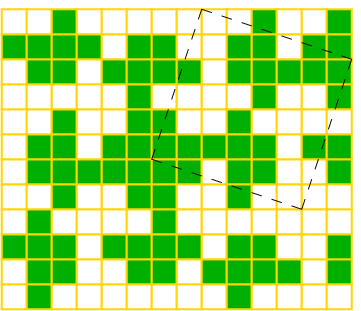}\hskip 10pt \includegraphics{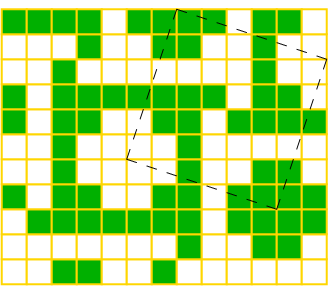}

(d) 4147-1*/787-1* \hskip 20pt (e) 4249-1*/2329-1*\hskip 20pt (f) 4249-2*/2329-2*

\vspace {4 pt}
\includegraphics{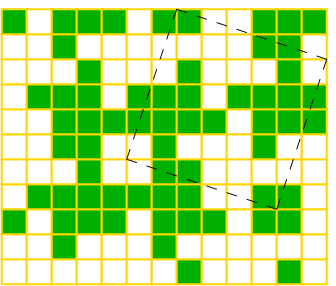}\hskip 10pt \includegraphics{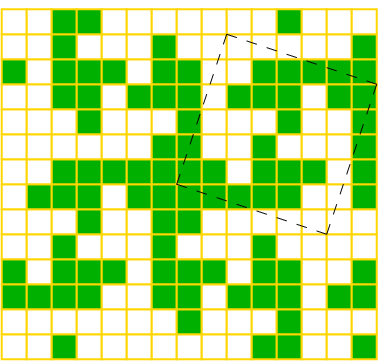}\hskip 10pt \includegraphics{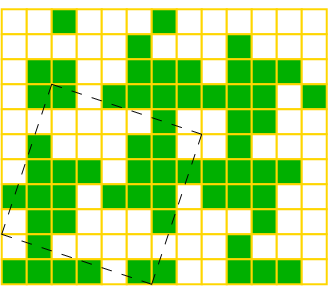}

(g) 4371-1*/4371-2*\hskip 25pt (h) 4371-2*/4371-1*\hskip 25pt (i) 4377-1*/4377-2*

\vspace {4 pt}
\includegraphics{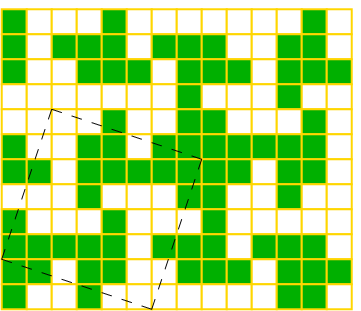}\hskip 10pt \includegraphics{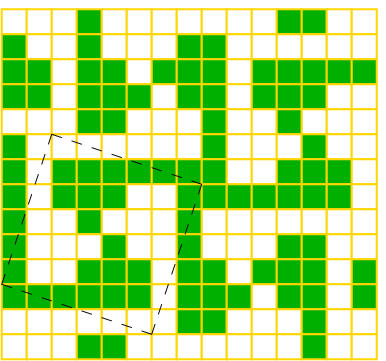}\hskip 10pt \includegraphics{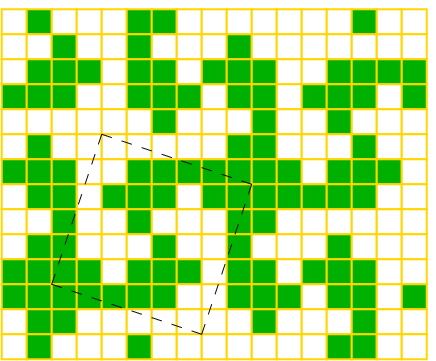}

(j) 4377-2*/4377-1*\hskip 30pt (k) 4387-1*/8367-1*\hskip 30pt (l) 4387-2*/8367-2*

\caption{Together with Figure 6 (787-2*, 4147-2*), the beginning of a catalogue of the order-20 species-37 isonemal prefabrics that fall apart, completed in Figures 11 and 12. The index number following each solidus is that of the reverse of the prefabric. The dashed squares are lattice units further explained in Section 6.}\label{fig10:}
\end{figure}
\begin{figure}
\centering

\includegraphics{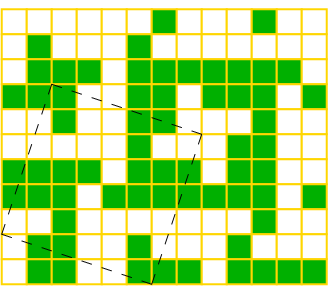}\hskip 10pt \includegraphics{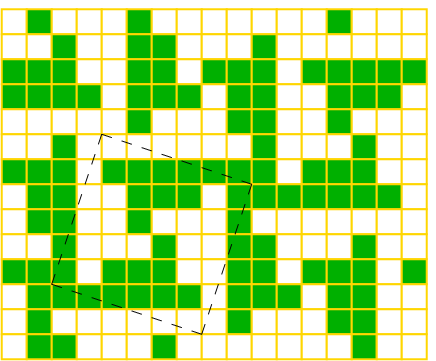}\hskip 10pt \includegraphics{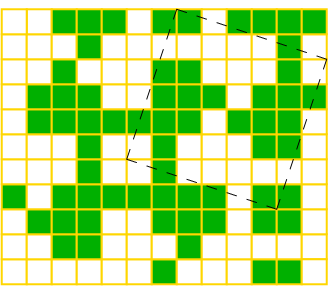}

(a) 4401-1*/4401-2*\hskip 25pt (b) 4401-2*/4401-1*\hskip 25pt (c)4489-1*/8497-1*

\vspace {4 pt}
\includegraphics{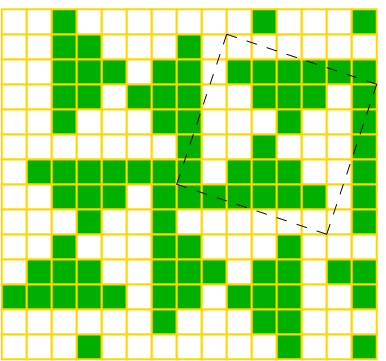}\hskip 10pt \includegraphics{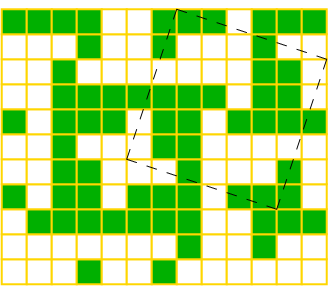}\hskip 10pt \includegraphics{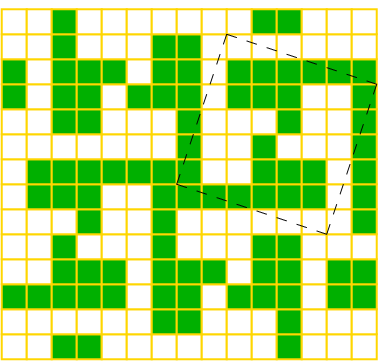}

(d) 4489-2*/8497-2*\hskip 25pt (e) 4643-1*/8723-1*\hskip 25pt (f) 4643-2*/8723-2*

\vspace {4 pt}
\includegraphics{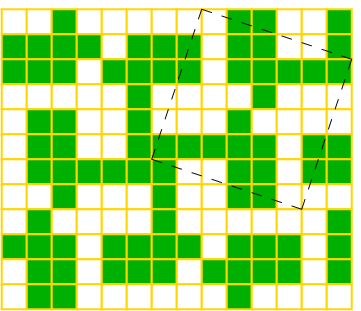}\hskip 10pt \includegraphics{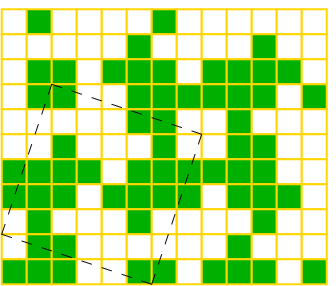}\hskip 10pt \includegraphics{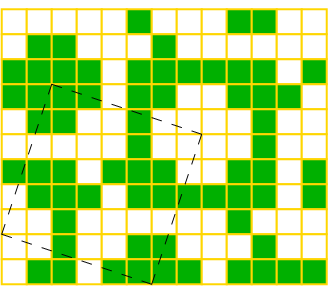}

(g) 8367-1*/4387-1*\hskip 20pt (h) 8367-2*/4387-2*\hskip 20pt (i) 8497-1*/4489-1*

\vspace {4 pt}
\includegraphics{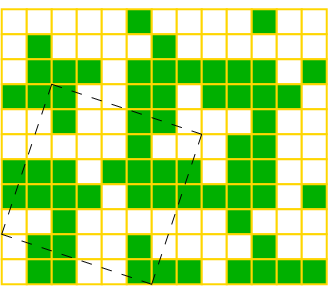}\hskip 10pt \includegraphics{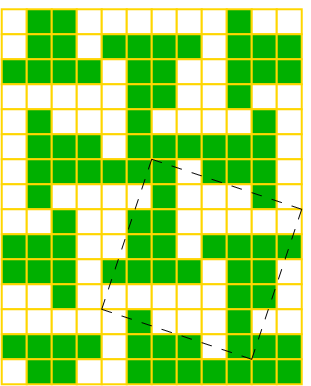}\hskip 10pt \includegraphics{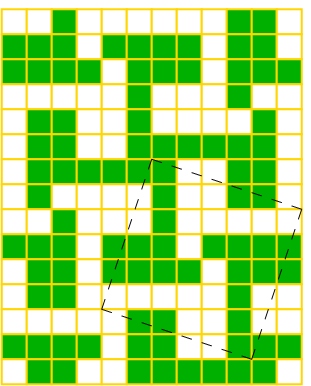}

(j) 8497-2*/4489-2*\hskip 20pt (k) 8723-1*/4643-1*\hskip 20pt (l) 8723-2*/4643-2*

\caption{Continuation of a catalogue of the order-20 species-37 isonemal prefabrics that fall apart.}\label{fig11:}
\end{figure}
\begin{figure}
\centering
\includegraphics{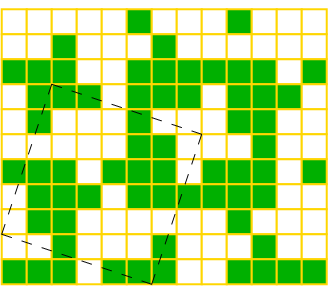}\hskip 10pt \includegraphics{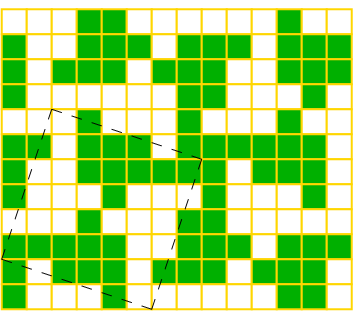}

(a) 8739-1*/8739-2*\hskip 20pt (b) 8739-2*/8739-1*

\vspace {4 pt}
\includegraphics{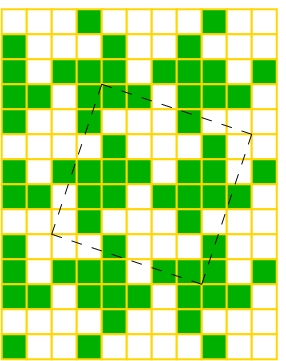}\hskip 20pt\includegraphics{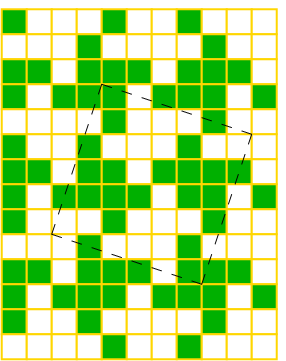}

(c) 34953-1*/34953-2*\hskip10pt (d) 34953-2*/34953-1*
\caption{Completion of a catalogue of the order-20 species-37 isonemal prefabrics that fall apart.}\label{fig12:}
\end{figure}

Figures 6 and 10--12 extend the catalogue of isonemal prefabrics that fall apart {\it for genus V} to order 20 and, since thick striping without quarter-turn symmetry occurs only in orders a multiple of 8 (Corollary 3.5), complete the extension (in \cite{P4}) for all prefabrics that fall apart to order 20.

One wonders naturally whether patterns of isonemal prefabrics with rotational symmetry that fall apart can all be produced as perfect colourings of isonemal fabrics, as has been shown not always to be possible in the case of prefabrics with axes of symmetry.
The exceptional prefabric 4-1-2* can be produced by striping the strands (though not perfectly) of 4-1-1, which is also exceptional, and in many other ways, one of which was described in the second-last paragraph.
The natural place to look for fabrics is at designs with the same lattice units and whose lattice units are preserved in the striping, species $33_4$.
We can show that every pattern that is the design of a species-37 prefabric that falls apart can be formed by striping a species-$33_4$ fabric.
\begin{thm}{Every non-exceptional design of an isonemal prefabric that has quarter-turn symmetry and that falls apart is the pattern produced by thickly striping the strands of a fabric of species $33_4$ or having symmetry group of Roth type $33_4$ as a subgroup.}
\end{thm} 

\begin{proof} This proof is long enough that the reader may benefit from seeing in advance how it goes. After some preliminaries to establish what is being discussed and how (para.~2, 3), the process for producing a fabric design from the design of a prefabric that falls apart is set out, first for the easier class of cells (para.~4) and then for the more interesting second half (para.~5--7). An example will be given of how the process works (para.~8). The symmetry group is established (para.~9) and then it is proved that the group does act as it must on the second half of the cells (para.~10) and then on the first half of them (para.~11--14).

Let there be the design of a non-exeptional isonemal genus-V prefabric that falls apart.
It must be of species 37.
Its appearance is the pattern from which the process must begin.
Pairs of rows and of columns are alternately predominantly (3/4) dark and pale. Where pale and pale and where dark and dark intersect are blocks that will be called, as though they had resulted from colouring, `redundant' by a handy slight abuse of term. 
Likewise the other half of the blocks will be called `irredundant'. 
There can be no doubt what apparent colouring is to be `reversed'.

There are four kinds of cell in the species-37 pattern that need to be sent back to a species-$33_4$ design.
They appear in blocks of four, redundant or irredundant and having a symmetry centre of one of two kinds in the centre or not.
This looks like six possibilities but is only four.
A \blbox in the pattern cannot lie in the centre of an irredundant block because it would be relating a predominantly dark row and predominantly pale column or vice versa.
The \blboxx s in the pattern lie therefore only in the centres of what must be redundant blocks in the colouring of the design. 
A \whbox in the pattern cannot lie in the centre of a redundant block because it cannot act on such a block.
The \whboxx s in the pattern lie therefore only in the centres of what must be irredundant blocks in the colouring of the design.
There are therefore only four kinds of block, redundant with and without \blbox and irredundant with and without \whboxx , the redundant and irredundant blocks forming a checkerboard pattern.

We now imagine a partial design (like Figure 3b) consisting of the irredundant blocks, duplicating the apparent colouring in the predominantly pale rows and reversing the apparent colouring (as with $\tau$) in the predominantly dark rows.
Irredundant blocks, with and without \whboxx , are now coloured as they must be in the design.
We shall have to see eventually that this colouring is consistent with a symmetry group of Roth type $33_4$.

The \blboxx s in the pattern, which lie in the centres of redundant blocks, must be converted into \whboxx s in the design.
For this to be possible, the blocks must be converted to blocks like a quarter of Figure 1a or its colour complement (i.e., like a $2\times 2$ matrix with a diagonal dark) in the design, because each is rotated within itself in the design by the new \whbox in its centre, and all of them the same, because they are rotated to one another by the \whboxx s of the pattern and design.
These two block designs are the only kinds invariant under \whbox in the centre.
When either is rotated a quarter turn in the same location or elsewhere by a new or old \whbox (no $\tau$), it is changed to the other by the rotation but back by the warp-weft colouring convention.
So each is a satisfactory way of weaving the blocks in the design accommodating new \whboxx s.
They need not be woven the same way as the irredundant blocks containing old \whboxx s, since the symmetry groups involved do not mix the two; they do not in the pattern and, being in the same locations, cannot in the design.
In fact, there is some freedom, but what is needed here is to show that something is possible.
Weaving them oppositely to the irredundant blocks containing old \whboxx s is an obvious way to prevent too much symmetry.
But the blocks containing the new \whboxx s can be supposed to be woven appropriately in either of those two ways.

Now we turn to the redundant blocks without \blboxx .
They cannot be coloured all pale, which would be the simple choice.%
\footnote{This was proposed in the inadequate algorithm of \cite{JA}, which happens to work on 4-1-2* because there are no such blocks in it.}
If one is coloured pale, then half of its images under the symmetry group will be dark and half pale because of the \whboxx .
There are many ways to weave the redundant blocks.
They need not be woven the same way as the redundant blocks formerly containing \blboxx s, since these too are not mixed.
Again there is some freedom.
Weaving them oppositely to redundant blocks containing \blboxx s is an obvious way to prevent too much symmetry (e.g., the production of plain weave).

There is no constraint at all on how to weave the single such block in the order-20 designs that have been used as examples (Figures 6 and 10--12).
(In larger-order designs, there is more freedom with more such blocks.
The next larger case has order 52 and three such blocks.)
Once that one block (or in larger examples those blocks) in a lattice unit is determined, however, the remainder are all images of it (or them) under the symmetry group.

The procedure described above has been carried out on the design of prefabric 20-8367-$2^*$ illustrated in Figure 11h. 
The result is shown in Figure 13.
\begin{figure}
\centering
\includegraphics{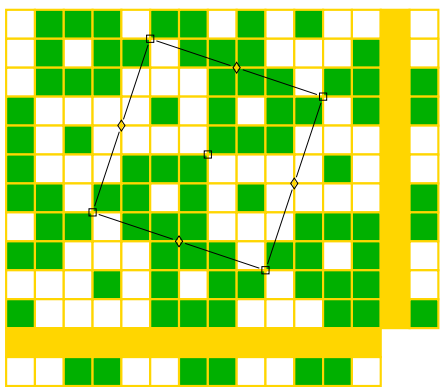}

\caption{Species-$33_4$ fabric with indication of thick striping to produce the design of the prefabric 20-8367-2* of Figure 11h.}\label{fig13:}
\end{figure}
The choices were made to weave redundant blocks containing new \whboxx s  (corners of illustrated lattice unit) the same way as the irredundant blocks containing the old \whboxx s and those not containing centres of rotation half dark and half pale with a straight boundary (where choice was free) for the sake of making an attractive coherent motif. 
The lattice unit marked in Figure 11h is marked in Figure 13 along with centres of rotation.

Type-$33_4$ symmetry operations in the design are determined by the \whboxx s of the pattern and those from converted \blboxx s in the pattern, and if the specified weave were coloured by thick striping the pattern resulting would be that of the given species-37 prefabric that falls apart.
What remains to be proved is that the weave specified must have the symmetry operations specified.

In the pattern, redundant blocks containing a \blbox are rotated onto redundant blocks containing a \blbox by both the \whboxx s and the \blboxx s, with a colour change in the former case and without a colour change in the latter.
In the design, these blocks are all woven the same way, and so they can be transformed among themselves by the \whboxx s in any way at all with no possible loss of consistency.
In the pattern, redundant blocks not containing a \blbox are transformed into redundant blocks not containing a \blbox by both the \whboxx s and the \blboxx s, with a colour change in the former case and without a colour change in the latter.
In the design, these blocks are all woven the same way, although probably not the same way as the other redundant blocks, and so they can be transformed among themselves consistently by the \whboxx s in any way at all.
All blocks that are redundant in the striping are woven consistently with the intended type-$33_4$ symmetry.

We need to see that the same is true of the irredundant blocks, the colouring of which was not arbitrarily determined.
The blocks that are irredundant in the pattern lie either in rows that are coloured predominantly dark in the striping (hence reversed in colour in the design) or those that are predominantly pale (and unchanged in colour in the design).
The \whboxx s in the design that were \whboxx s in the pattern, lying there in irredundant blocks, are either in pale rows and dark columns or dark rows and pale columns.
The \whboxx s in the design that were \blboxx s in the pattern, lying there in redundant blocks, are either in pale rows and pale columns or dark rows and dark columns.
So we have two kinds of block to show are treated consistently by old \whboxx s and new \whboxx s in two different kinds of position each, eight things to see.
These four positions of \whboxx s and \blboxx s  are set out as the rows of Table 2, where the two kinds of irredundant block are assigned in the next paragraph to the two columns of entries.

Let the co-ordinates of a centre of rotation be (0, 0).
Index the block with centre $2i$ cells to the right and $2j$ cells up ($i, j$, integers, not necessarily positive) by $(x, y) = (2i, 2j)$.
We shall be concerned with the parities of $i$ and $j$, $x$ and $y$ always being even.
If the origin is in an irredundant block, then the block is irredundant if the parity of $i$ and $j$ are the same, redundant if different.
If the origin is in a redundant block, then the block is redundant if the parity of the $i$ and $j$ are the same, irredundant if different.
Under a quarter turn, each block $(2i, 2j)$ is rotated to a block $(-2j, 2i)$, equality or inequality of parities of $i$ and $j$ still the same, therefore irredundancy or redundancy is preserved as we knew but see confirmed.
To consider the eight situations we face, Table 2 sets out the various parity combinations that occur in those situations.
\begin{table}
\caption{Parities of $i$ and $j$ of irredundant blocks for various centres of rotation.}\label{tab:2}
\centering
\begin{tabular}{lcc}
\toprule
&\multicolumn{2}{c}{Location of irredundant block}\\
\cmidrule(l){2-3}
Centre type & Dark row, pale column & Pale row, dark column\\
and location &Colour changed in design & Colour unchanged in design\\
\midrule
\whbox from \whbox in pattern\\
Pale row, dark column & $i, j$ odd & $i, j$ even \\ \addlinespace[2pt]
\whbox from \whbox in pattern\\
Dark row, pale column & $i, j$ even & $i, j$ odd \\ \addlinespace[2pt]
\whbox from \blbox in pattern\\
Pale row, pale column & $i$ even, $j$ odd & $i$ odd, $j$ even\\ \addlinespace[2pt]
\whbox from \blbox in pattern\\
Dark row, dark column & $i$ odd, $j$ even & $i$ even,$j$ odd \\
\bottomrule
\end{tabular}
\end{table}

In the top two rows of entries, where the centre is in an irredundant block, while the parities are different in different positions, those of $i$ and $j$ are the same.
We see that the result of a quarter-turn (includes reversing $x$ and $y$ indices) is to take changed-colour blocks to changed-colour blocks and unchanged-colour blocks to unchanged-colour blocks because the parities are preserved; no column change in Table 2.
The action of a \whbox in the design arising from a \whbox in the pattern is the same as the action of the \whbox in the pattern.
If the image was changed in the design, then so was the pre-image, and if not then not.
The \whboxx s arising from \whboxx s in the pattern are indeed symmetries of the design.

In the bottom two rows of the table, where the centre is in a redundant block, the parities differ in all positions.
The result of a quarter-turn (reversing $x$ and $y$ indices) is to take changed-colour blocks to unchanged-colour blocks and unchanged-colour blocks to changed-colour blocks because the positions of odd and even $i$ and $j$ are reversed; columns in Table 2 are swapped.
The action of a \whbox in the design arising from a \blbox in the pattern is that of the \blbox in the pattern (\blbox being the composition of \whbox and $\tau$) composed with a second $\tau$.
The \whboxx s of the design arising from \blboxx s of the pattern are symmetries of the design because the $\tau$s cancel each other.

Like the redundant blocks, the irredundant blocks are transformed as by a group of Roth type $33_4$. The design is that of a fabric of species $33_4$ or has a group of Roth type $33_4$ as a subgroup.
\end{proof}

If the design produced is not of species $33_4$, then the thick striping will not be a perfect colouring.
If, however, a fabric with `too much' symmetry can be produced, symmetries can probably be broken to reduce the symmetry group to type $33_4$.
\section{Woven Cubes}
\noindent This section illustrates that the thick striping of strands is relevant to what I call the perfect colouring of woven cubes. It was shown in \cite{P3} that lattice units of an isonemal prefabric of order greater than 4 could be the regions of the net of an isonemal woven cube if and only if the prefabric is of species 33, 34, 37, 38, or 39.
Woven cubes were introduced into the weaving literature by Jean Pedersen (e.g., \cite{JP}).
The focus in \cite{JP, P3} is on the isonemality of the (essentially colourless) weaving, not on what such a cube would look like.
In \cite{P3} the normal colouring of prefabric designs was used in the usual way to represent nets of cubes, but it makes no sense to think of a cube's being woven from `normally' coloured strands since the very idea of warp and weft as distinct strands makes no sense in the context.
In Shephard's simplest case with oblique lattice units \cite[Figure 17]{JP} and so having faces of area 5, each of the 6 strands required to weave the cube makes an appearance on all 6 faces; even there separation into warps and wefts makes no sense.
In examples even slightly larger, strands cross {\it themselves} (in the cube of Figures 14 and 15, each of the 8 strands crosses itself 3 times).
While it may be possible to colour strands other than by striping warp and weft in the net of the cube, that is the simplest thing to try, and it can work.
\begin{figure}
\centering
\includegraphics{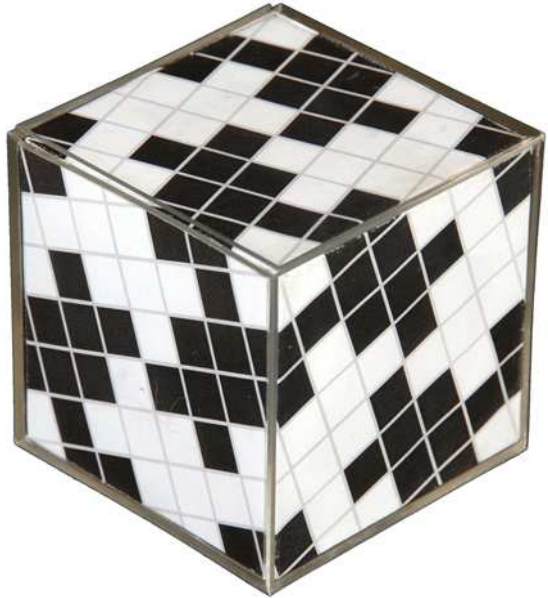}

\caption{Downward view of cube coloured by the thick striping catalogued as 20-8367-2* showing three faces.}\label{fig14:}
\end{figure}

The definition of isonemal woven cube used in \cite{P3} requires not only that the six faces be woven in the same way such that the rotations of the cube around its 4-fold, 3-fold, and 2-fold axes be collectively transitive on the strands but also that there be no centres of symmetry in the net of the cube other than the quarter-turns at the centres and corners of faces and the half-turns at the middle of their sides. 
Obviously it is possible to weave cubes from prefabrics with more symmetries than can be given by the rotations of the whole cube (e.g., plain weave in \cite[Figure 17]{JP}).

Pedersen expressed reservations about applying reflection in the plane of the fabric to operations on the net of the cube, but in \cite[Figure 17]{JP} and \cite{P3} those reservations were ignored.
In this paper, however, we have found that the presence of quarter-turns requiring $\tau$ leads to appearances of obverse and reverse with contrasting stripes so related.
As a result, while it is possible to colour perfectly fabrics of species $36_s$ thinly and $33_4$, $35_4$, and 37 thickly, we know that the results in the cases of $36_s$ and $35_4$ do not look symmetrical, and they are also not ways of weaving a cube.

The species to examine are $33_4$ and 37.

\cite[Lemma 10]{P3} shows that the isonemality of a woven cube requires that there be no \blbox in the net at face corners.
That means that, if a species-37 prefabric is to be used to weave a cube, the net of the cube has corners at the \whbox of Figures 8a and 9a rather than being composed of lattice units marked in those figures, putting \blbox in the block (irredundant in Figure 8a, redundant in Figure 9a) in the centre of every face.
Note that the colouring convention does not operate in these figures; they are patterns not designs.
Each face accordingly has a stripe across it; its quarter-turn symmetry is not a colour symmetry.
Moreover, in Figure 9a, where the face-centre blocks are redundant, the four cells surrounding the corner \whbox are differently coloured, two dark and two pale, so that the net must include two of one colour and one of the other, destroying in the cube the perfect symmetry in the plane.

In contrast, thick striping of species $33_4$ can be used to colour cubes provided only that the centres of redundant blocks are used as the corners of the cube, putting \whbox in the centre of irredundant blocks in the centre of faces.
These faces therefore look like the lattice units illustrated in Figures 6 and 10--12.
All the faces look the same.

\begin{defn}
A {\em perfectly coloured isonemal woven cube} is a cube woven isonemally \cite{P3} with its strands so coloured that each rotational symmetry of the cube preserves or permutes the colours showing in the cells.
\end{defn}
\begin{thm}
Lattice units of a fabric of species $33_4$, perfectly coloured by thick striping, compose the net of a perfectly coloured isonemal woven cube provided that the corners of the lattice units fall at the centres of redundant blocks of the striping.
\end{thm}
\begin{proof}
We know that thick striping must have centres of redundant blocks at either corners or centres of lattice units.
For the lattice units to compose the net of a perfectly coloured cube, the centres of redundant blocks must be placed at corners so that, of the four cells in a block, the three that fall in the net have the same colour; if a corner block were irredundant, containing as it does a centre of quarter-turn, two of its cells would be dark and two pale.
We know that the thick striping specified is a perfect colouring of the planar fabric and that the weaving makes an isonemal cube.
What needs to be shown is that the rotations of the cube permute or preserve the colours showing in the cells.
The quarter-turns at the centres of the faces obviously reverse the colours because the centre blocks are irredundant and woven accordingly so that quarter-turns reversed colours in the plane.
The half-turns at the middle of edges obviously reverse the colours because they do so in the plane, lying as they do on the strand boundary between thick stripes predominantly dark and predominantly pale.
The rotations whose working is not obvious are the quarter-turns at the corners of the lattice units.
They are quarter-turns without $\tau$ in order to make the weaving isonemal, but their blocks are made redundant in order to make their effect look like quarter-turns with $\tau$ in standardly coloured design diagrams, i.e., not to reverse colours.
Because colours are preserved in the quarter-turn at a lattice-unit corner that is a vertex $P$ of the cube, predominantly dark pairs of strands are rotated to predominantly dark pairs of strands and likewise pale {\it in the plane}.
In all four directions, predominantly dark pairs of strands are the same distances from $P$, likewise pale pairs.
In the net, where one of the four lattice units surrounding $P$ is removed, when the formerly planar strands are joined up along lattice-unit boundaries to make a new vertex $Q$ of the cube and the edge $PQ$, predominantly dark pairs match with predominantly dark pairs and likewise pale because they are the same distance from $P$.
The result of the joining along $PQ$ looks exactly the same as the result of mere folding along the other edges, say, $PR$ and $PS$.

The centre of colour-preserving 4-fold rotation at $P$ in the plane now lies on an axis of colour-preserving 3-fold rotation of the cube in space.
$Q$, $R$, and $S$ fall at the centres of redundant blocks coloured complementary to $P$'s block because the edge lengths arrange it so.
Perfect colouring of the cube is the case.
\end{proof}

The 4-sided picture-frame-like motifs of 4147-2* (Figure 6b), 2329-1*, 2*, 4147-1*, and 4249-2* (in Figure 10), 8367-2* (Figure 11h --- also 14 and 15), and a number of others, which surround such centres as $P$ of the theorem become 3-sided in the cube, which looks odd (Escher-like) in two-dimensional representations since they are still composed of square cells.

Figures 14 and 15 illustrate the thick striping catalogued as 20-8367-2* used as the net of a cube. 
To have a fabric to stripe, it is necessary to choose one of the ways to weave a species-$33_4$ fabric by reversing the colouring algorithm and fixing the weaving of the redundant blocks compatibly with the symmetry group.
This was done in Figure 13.
The lattice unit marked in Figures 11h and 13 becomes the top surface of the cube in Figure 14.
Figure 14 shows three faces of the cube with a picture-frame motif around a vertex at the front in what is left of a pale redundant block inside the picture frame.
Because the vertex at the intersection of the three hidden faces is a dark vertex, the faces look like the colour-complement of Figure 14 upside down.
For the same reason that reverses of fabrics are shown reflected, namely the matching of corresponding cells, Figure 15 illustrates the hidden faces reflected in a mirror.
\begin{figure}
\centering
\includegraphics{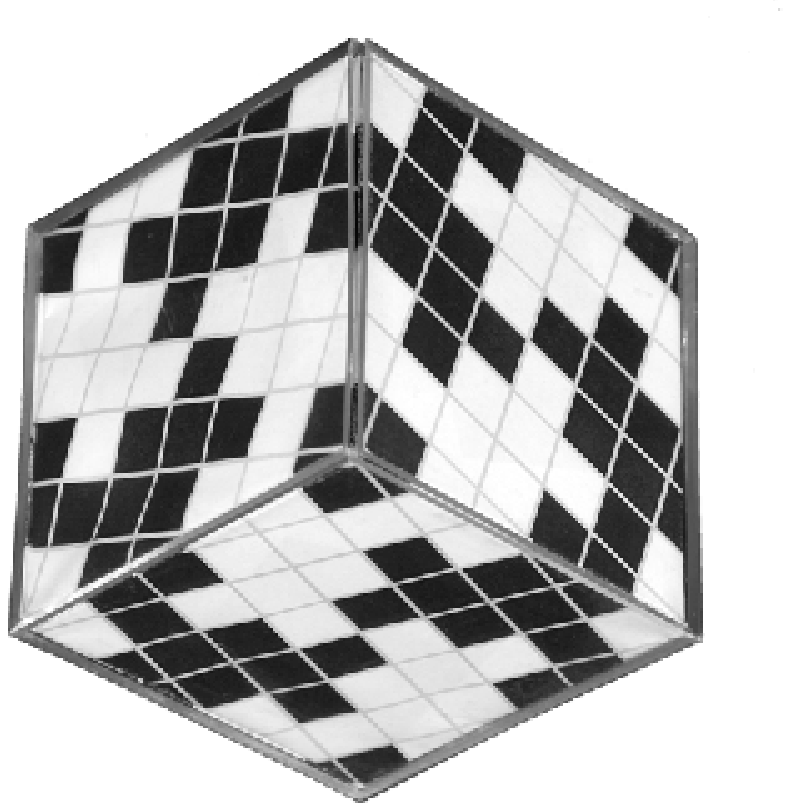}

\caption{Upward view of the back of the cube of Figure 13 reflected in a mirror.}\label{fig15:}
\end{figure}
Reflection has the effect that the cell boundaries around the edges of Figures 14 and 15 match, allowing one to follow a strand off the edge of the illustration in either figure to where it continues (within the same cell) in the other figure.
For example, the dark strand that can be seen at the bottom front corner of Figure 14 proceeding upward and to the right, is hidden for a cell then exposed for a cell, then passes through a dark redundant block before being hidden for one cell and reaching an exposed cell that overlaps the edge of the cube.
It can be seen again, as the third strand up the right vertical edge in Figure 15, to proceed up and to the left through another redundant block, exposed for 1 cell, hidden for 1, through 2 redundant blocks and exposed in the 2 cells between them as it passes across the central vertical edge, and then hidden for 2 cells before reaching the upper left vertex as a cell in a redundant block mainly invisible in Figure 15.
The strand is transferred to Figure 14 with the second (larger) part of that cell in the larger part of the same redundant block and a small part of the next cell on the top face coming forward before descending the lower left face.
In this way it is possible to follow the whole course of the strand---or any strand---as it makes its way twice in perpendicular directions across each face.
One sees that the strand is 60 cells long (3 times its 1-dimensional period), that there are 4 dark and 4 pale strands, their 480 cells being a double covering of the 240 cells of the cube.

\ack
{Work on this material has been done at home and at Wolfson College, Oxford.
Richard Roth helped with the understanding of his papers.
Will Gibson made it possible for me to draw the diagrams with surprising ease from exclusively keyboard input. 
Allen Patterson of Information Services and Technology at the University of Manitoba spent a lot of time on the photographs.
To them and Wolfson College I make grateful acknowledgement.}


\end{document}